\documentclass[12pt]{amsproc}
\usepackage[cp1251]{inputenc}
\usepackage[T2A]{fontenc} 
\usepackage[russian]{babel} 
\usepackage{amsbsy,amsmath,amssymb,amsthm,graphicx,color} 
\usepackage{amscd}
\def\le{\leqslant}
\def\ge{\geqslant}

\newtheorem{prop}{Предложение}

\theoremstyle{definition}

\theoremstyle{remark}

\begin {document}
\unitlength=1mm
\title[Совместные формулы умножения]
{Совместные формулы умножения для многочленов Апостола-Бернулли и обобщённых многочленов Фробениуса-Эйлера}
\author{Г. Г. Ильюта}
\email{ilgena@rambler.ru}
\address{}
\thanks{Работа поддержана грантом РФФИ-20-01-00579}

%\centerline{Содержание}

\begin{abstract}
Мы докажем совместные формулы умножения для многочленов Апостола-Бернулли и обобщённых многочленов Фробениуса-Эйлера и свяжем их с суммами Дедекинда-Радемахера, суммами Апостола-Дедекинда и суммами Фурье-Дедекинда.

We prove the simultaneous multiplication formulas for Apostol-Bernoulli polynomials and generalized Frobenius-Euler polynomials and relate them to Dedekind-Rademacher sums, Apostol-Dedekind sums and Fourier–Dedekind sums.
\end{abstract}

\maketitle
\tableofcontents

\bigskip

\section{Введение}

\bigskip

  Мы докажем совместные формулы умножения для двух примеров многочленов Аппеля -- многочленов Апостола-Бернулли $B_m(q,\lambda)$ и обобщённых многочленов Фробениуса-Эйлера $H_m^{(p)}(q,\lambda,\gamma)$, которые определяются производящими функциями 
$$
\frac{t}{\lambda e^t-1}e^{qt}=\sum_{m=0}^{\infty}B_m(q,\lambda)\frac{t^m}{m!},
$$
$$
\frac{(1-\gamma)^p}{\lambda e^t-\gamma}e^{qt}=\sum_{m=0}^{\infty}H_m^{(p)}(q,\lambda,\gamma)\frac{t^m}{m!},  
$$
где $p\in\mathbb Z$. Мы интерпретируем эти формулы как соотношения между суммами Дедекинда, так как обобщённые многочлены Фробениуса-Эйлера появляются в них через суммы Дедекинда
$$
E_{m,n}^{r,p}(q,\lambda;C):=\sum_{k=1}^{n-1}\frac{\epsilon_n^{-kr}H_{m-1}^{(p)}(q,\lambda,\epsilon_n^{-k})C_{-k}}{(1-\epsilon_n^k)^p},						\eqno (1)
$$
где $m\in\mathbb Z_{>0}$, $n\in\mathbb Z_{>1}$, $r,p\in\mathbb Z$, $\epsilon_n:=e^{\frac{2\pi i}{n}}$, $C=\{C_k\}$ -- $n$-периодическая последовательность. Такая терминология объясняется следующими двумя примерами: при $m=p=\lambda=1$, 
$$
C_k=\frac{1}{1-\epsilon_n^{-ak}},\quad 0<k<n,
$$
в определении (1) приходим к представлению сумм Дедекинда-Радемахера в виде сумм Фурье-Дедекинда \cite{1}, p.~156, а  при $m>1$, $r=q=0$, $p=\lambda=1$,
$$
C_k=\frac{1}{1-\epsilon_n^{ak}},\quad 0<k<n,
$$
 -- к аналогичному представлению для сумм Апостола-Дедекинда \cite{2}, p.~521.

\bigskip

\section{Дискретное преобразование Фурье}

\bigskip

 Пусть $n$-периодическая последовательность $K=\{K_j\}$ связана с последовательностью $C=\{C_k\}$ дискретным преобразованием Фурье
$$
C_k=\sum_{j=0}^{n-1}K_j\epsilon_n^{kj},
$$
$$
K_j=\frac{1}{n}\sum_{k=0}^{n-1}C_k\epsilon_n^{-kj}.
$$

  Интерполяционный многочлен Лагранжа $C^{(r)}(q)$ последовательности $\{\epsilon_n^{-kr}C_{-k}\}$ определим условиями $\deg C^{(r)}(q)<n$ и 
$$
C^{(r)}(\epsilon_n^{-k})=\epsilon_n^{-kr}C_{-k},\quad 0\le k\le n-1.
$$

\begin{prop}\label{prop1} 
$$
C^{(r)}(q)=\sum_{j=0}^{n-1}K_{j-r}q^j.				\eqno (2)	
$$
\end{prop}

  Доказательство. По формуле Лагранжа
$$
\frac{C^{(r)}(q)}{q^n-1}=\frac{1}{n}\sum_{k=0}^{n-1}\frac{\epsilon_n^{-k}\epsilon_n^{-kr}C_{-k}}{q-\epsilon_n^{-k}}						\eqno (3)
$$
$$
=\frac{1}{n}\sum_{k=0}^{n-1}\frac{\epsilon_n^{-k}\sum_{j=0}^{n-1}K_{j-r}\epsilon_n^{-kj}}{q-\epsilon_n^{-k}}
$$
$$
=\frac{1}{n}\sum_{j=0}^{n-1}K_{j-r}\sum_{k=0}^{n-1}\frac{\epsilon_n^{-k}\epsilon_n^{-kj}}{q-\epsilon_n^{-k}}.
$$
$$
=\frac{\sum_{j=0}^{n-1}K_{j-r}q^j}{q^n-1}.\qquad\Box
$$

\bigskip

\section{Совместные формулы умножения}

\bigskip

\begin{prop}\label{prop2} Для $m\in \mathbb Z_{\ge 1}$
$$
(-1)^{p-1}mE_{m,n}^{r,p}(nq,\lambda;C)
$$
$$
=C_0B_m(nq,\lambda)-n^m\sum_{j=0}^{n-1}K_{j-r-p+1}\lambda^jB_m\left(q+\frac{j}{n},\lambda^n\right).								\eqno (4)	
$$
\end{prop}
 
   Доказательство. Наличие производящей функции для обобщённых многочленов Фробениуса-Эйлера позволяет написать аналогичную формулу для сумм (1)
$$
G_n^{r,p}(q,\lambda,t;C):=\sum_{m=0}^{\infty}E_{m+1,n}^{r,p}(q,\lambda;C)\frac{t^m}{m!}
$$
$$
=\sum_{m=0}^{\infty}\sum_{k=1}^{n-1}\frac{\epsilon_n^{-kr}H_m^{(p)}(q,\lambda,\epsilon_n^{-k})C_{-k}}{(1-\epsilon_n^k)^p}\frac{t^m}{m!}
$$
$$
=\sum_{k=1}^{n-1}\frac{\epsilon_n^{-kr}C_{-k}}{(1-\epsilon_n^k)^p}\sum_{m=0}^{\infty}H_m^{(p)}(q,\lambda,\epsilon_n^{-k})\frac{t^m}{m!}
$$
$$
=\sum_{k=1}^{n-1}\frac{\epsilon_n^{-kr}C_{-k}(1-\epsilon_n^{-k})^p}{(1-\epsilon_n^k)^p(\lambda e^t-\epsilon_n^{-k})}e^{qt}
$$
$$
=(-1)^p\sum_{k=1}^{n-1}\frac{\epsilon_n^{-k}\epsilon_n^{-k(r+p-1)}C_{-k}}{(\lambda e^t-\epsilon_n^{-k})}e^{qt}.							\eqno (5)
$$
Добавляя в последнюю сумму слагаемое, отвечающее $k=0$, и сравнивая её с формулой (3), получим
$$
\frac{C_0}{\lambda e^t-1}e^{qt}+(-1)^pG_n^{r,p}(q,\lambda,t;C)=\frac{nC^{(r+p-1)}(\lambda e^t)}{\lambda^ne^{nt}-1}e^{qt}.
$$
Умножаем это соотношение на $te^{(n-1)qt}$ и применяем формулу (2)
$$
\frac{C_0t}{\lambda e^t-1}e^{nqt}+(-1)^ptG_n^{r,p}(q,\lambda,t;C)e^{(n-1)qt}
$$
$$
=\frac{nt\sum_{j=0}^{n-1}K_{j-r-p+1}\lambda^je^{jt}}{\lambda^ne^{nt}-1}e^{nqt}=\frac{nt\sum_{j=0}^{n-1}K_{j-r-p+1}\lambda^je^{\left(q+\frac{j}{n}\right)nt}}{\lambda^ne^{nt}-1}.  \eqno (6)
$$
Проходя цепочку равенств (5) в обратном направлении, имеем
$$
tG_n^{r,p}(q,\lambda,t;C)e^{(n-1)qt}=t(-1)^p\sum_{k=1}^{n-1}\frac{\epsilon_n^{-k}\epsilon_n^{-k(r+p-1)}C_{-k}}{(\lambda e^t-\epsilon_n^{-k})}e^{nqt}
$$
$$
=t\sum_{m=0}^{\infty}E_{m+1,n}^{r,p}(nq,\lambda;C)\frac{t^m}{m!}=\sum_{m=1}^{\infty}mE_{m,n}^{r,p}(nq,\lambda;C)\frac{t^m}{m!}.						\eqno (7)
$$
Раскладывая в формуле (6) производящие функции для многочленов Апостола-Бернулли в ряды и применяя формулу (7), получим
$$
C_0\sum_{m=0}^{\infty}B_m(nq,\lambda)\frac{t^m}{m!}+(-1)^p\sum_{m=1}^{\infty}mE_{m,n}^{r,p}(nq,\lambda;C)\frac{t^m}{m!}
$$
$$
=\sum_{j=0}^{n-1}K_{j-r-p+1}\lambda^j\sum_{m=0}^{\infty}n^mB_m\left(q+\frac{j}{n},\lambda^n\right)\frac{t^m}{m!}
$$
и теперь сравниваем коэффициенты при одинаковых степенях переменной $t$.$\qquad\Box$

  Полагая в формуле (4) $C_0=n$, $C_k=0$, $n\nmid k$,
$$
K_j=\frac{1}{n}\sum_{k=0}^{n-1}C_k\epsilon_n^{-kj}=1, \quad j\in\mathbb Z,
$$
$$
E_{m,n}^{r,p}(q,\lambda;C)=\sum_{k=1}^{n-1}\frac{\epsilon_n^{-kr}H_{m-1}^{(p)}(q,\lambda,\epsilon_n^{-k})C_{-k}}{(1-\epsilon_n^k)^p}=0,
$$
приходим к формулам умножения для многочленов Апостола-Бернулли \cite{3}
$$
B_m(nq,\lambda)=n^{m-1}\sum_{j=0}^{n-1}\lambda^jB_m\left(q+\frac{j}{n},\lambda^n\right).	
$$

\bigskip

\section{Пример}

\bigskip

  Пусть для $1\le j\le n$ $K_j=1$, если $(j,n)=1$, и $K_j=0$, если $(j,n)>1$, где $n\in\mathbb Z_{>1}$ и $(j,n)$ -- наибольший общий делитель. Тогда
$$
C_k=\sum_{j=0}^{n-1}K_j\epsilon_n^{kj}=\sum_{\substack{1\le j\le n\\(j,n)=1}}\epsilon_n^{kj}
$$
-- суммы Рамануджана (степенные суммы корней циклотомического многочлена), в частности, $C_0=\phi(n)$ -- функция Эйлера, $C_1=\mu(n)$ -- функция Мёбиуса. Согласно \cite{5}, p. 75,
$$
C^{(0)}(q)=\sum_{j=0}^{n-1}K_{j}q^j
$$
$$
=\sum_{\substack{1\le j\le n\\(j,n)=1}}q^j=(q^n-1)\sum_{d|n}\frac{\mu(d)q^d}{q^d-1}=(q^n-1)\sum_{d|n}\frac{\mu(d)}{q^d-1}.
$$
Последнее равенство следует из условия $n>1$ и определения функции Мёбиуса: $\mu(1)=1$,
$$
\sum_{d|n}\mu(d)=0,\quad n>1.
$$
Пусть
$$
V_n^{(k)}(\lambda):=\sum_{\substack{1\le j\le n\\(j,n)=1}}j^k\lambda^j=\left(q\frac{d}{dq}\right)^k\sum_{j=0}^{n-1}K_{j}q^j.		
$$
Ссылки на источники, содержащие явные формулы для сумм $V_n^{(k)}(1)$, имеются в \cite{4}, p.~242. Например, $V_n^{(0)}(1)=\phi(n)$ и
$$
V_n^{(1)}(1)=\sum_{\substack{1\le j\le n\\(j,n)=1}}j=\frac{1}{2}\sum_{\substack{1\le j\le n\\(j,n)=1}}(j+(n-j))=\frac{n\phi(n)}{2}.		
$$

  При $r+p=1$ формула (4) примет вид
$$
(-1)^{p-1}mE_{m,n}^{r,p}(nq,\lambda;C)
$$
$$
=\phi(n)B_m(nq,\lambda)-n^m\sum_{\substack{1\le j\le n\\(j,n)=1}}\lambda^jB_m\left(q+\frac{j}{n},\lambda^n\right).					$$
Как известно
$$
B_m(q,\lambda)=\sum_{i=0}^m\binom{m}{i}q^{m-i}B_i(\lambda), 
$$
где $B_i(\lambda):=B_i(0,\lambda)$. Поэтому
$$
-\frac{(-1)^{p-1}mE_{m,n}^{r,p}(nq,\lambda;C)-\phi(n)B_m(nq,\lambda)}{n^m}=\sum_{\substack{1\le j\le n\\(j,n)=1}}\lambda^jB_m\left(q+\frac{j}{n},\lambda^n\right)					
$$
$$
=\sum_{\substack{1\le j\le n\\(j,n)=1}}\lambda^j\sum_{i=0}^m\binom{m}{i}\left(q+\frac{j}{n}\right)^{m-i}B_i(\lambda^n) 
$$
$$
=\sum_{i=0}^m\binom{m}{i}B_i(\lambda^n)\sum_{\substack{1\le j\le n\\(j,n)=1}}\lambda^j\sum_{k=0}^{m-i}\binom{m-i}{k}\frac{j^k}{n^k}q^{m-i-k}
$$
$$
=\sum_{i=0}^m\sum_{k=0}^{m-i}\frac{1}{n^k}\binom{m}{i}\binom{m-i}{k}B_i(\lambda^n)q^{m-i-k}\sum_{\substack{1\le j\le n\\(j,n)=1}}j^k\lambda^j
$$
$$
=\sum_{i=0}^m\sum_{k=0}^{m-i}\frac{m!}{n^ki!k!}B_i(\lambda^n)V_n^{(k)}(\lambda)\frac{q^{m-i-k}}{(m-i-k)!}.
$$
В результате
$$
(-1)^{p-1}mE_{m,n}^{r,p}(nq,\lambda;C)
$$
$$
=\phi(n)B_m(nq,\lambda)-\sum_{i=0}^m\sum_{k=0}^{m-i}\frac{n^{m-k}m!}{i!k!}B_i(\lambda^n)V_n^{(k)}(\lambda)\frac{q^{m-i-k}}{(m-i-k)!}.				
$$

\end {document}